%% file: Euler_residual_note.tex
\documentclass[reqno]{amsart}

\usepackage[top=1in, bottom = 1in, left = 1.16in, right = 1.16in]{geometry}
 \usepackage{amsmath,amssymb}
 \usepackage{algpseudocode}
 \usepackage{algorithm}
 \usepackage{algorithmicx}
 \usepackage{caption}
 \usepackage{subcaption}
 \usepackage{amsthm}
 \usepackage{esint}
 \numberwithin{equation}{section}
 \usepackage{amsfonts}
 \usepackage{graphicx}
 \usepackage{hyperref}
\usepackage{enumitem}
 \usepackage{makecell}
 \usepackage{longtable}
 \usepackage{forest}
 \usepackage{comment}
\input{JCmacro}

\def \lloc {\mathsf{loc}}
\def \nloc {\mathsf{nloc}}

\newcommand\jan[1]{  { \la{#1}\ra } }

\author{Jiajie Chen}
\thanks{Department of Mathematics, University of Chicago, Chicago, IL 60637. Email: jiajiechen@uchicago.edu}
\author{Thomas Y. Hou}
\thanks{Applied and Computational Mathematics, Caltech, Pasadena, CA 91125. Email:hou@cms.caltech.edu.}
\date{\today}

\title[Clarification on the distinction between two errors]{A clarification on the distinction between nonlocal error and profile residual error in a computer-assisted proof of 3D Euler singularity}

\begin{document}

\begin{abstract}
In \cite{zhang2025dimension}, the author discusses the applicability of a
nonexistence result for self-similar profiles to the Hou--Luo scenario and
raises questions about the approximate self-similar profile constructed in
\cite{ChenHou2023a,ChenHou2023b}. In this note, we clarify two main points. First,
the far-field discrepancy identified in \cite{zhang2025dimension} is the raw
nonlocal Poisson error, not the profile residual error used in the stability
proof of \cite{ChenHou2023a,ChenHou2023b}. The relevant nonlocal contribution
to the profile residual contains the nonlocal velocity error multiplied by
additional decaying factors, which provide crucial smallness in the far field estimates. 
Second, we examine the figure-based evidence used in
\cite[Remark 2.4]{zhang2025dimension} to support the applicability of the
assumption in \cite[Proposition 2.3(i)]{zhang2025dimension}, and show that the same grid-point data do not support that interpretation. Thus the comparisons and figure-based evidence in \cite{zhang2025dimension} neither invalidate the residual estimates in \cite{ChenHou2023a,ChenHou2023b} nor justify the claimed applicability of the angular-increase assumption to the Hou--Luo scenario.

\end{abstract}

\maketitle

\section{Introduction}

Self-similar methods play an important role in the study of finite time singularities 
in incompressible fluid equations \cite{ChenHou2023a,ChenHou2023b,elgindi2019finite,Tristan2022,chen2019finite2}. To construct blowup, an important approach is to first construct an approximate or exact blowup profile and then prove its nonlinear stability \cite{chen2019finite,chen2021HL,ChenHou2023a,ChenHou2023b}.

Recently, under various assumptions, Zhang \cite{zhang2025dimension} proved nonexistence results for certain self-similar profiles for 
axisymmetric Euler equations. He discussed the applicability of these nonexistence results to the Hou--Luo scenario \cite{luo2013potentially-2,ChenHou2023a,ChenHou2023b}. Moreover, \cite{zhang2025dimension} raised a question about a possible discrepancy involving the nonlocal error in \cite{ChenHou2023a,ChenHou2023b}. In this short note, we clarify that this comparison rests on an incorrect identification of the raw nonlocal Poisson error with the profile residual error used in the stability analysis. Once this distinction is made, the comparison in \cite{zhang2025dimension} does not provide a valid objection to, nor can it be used to cast doubt on, the residual estimates in \cite{ChenHou2023a,ChenHou2023b}. We do not discuss here the main dimension-reduction results of \cite{zhang2025dimension}; our purpose is only to clarify the specific numerical and interpretive issues concerning the Hou--Luo profile discussed in Remarks 2.4--2.6 of \cite{zhang2025dimension}.

\section{Profile residual error and nonlocal error }

Recall the dynamic rescaling equations of the Boussinesq equation in $\R^2_+$ 
from \cite[Section 2.1]{ChenHou2023a}
\beq\label{eq:bousdy1}
\bal
\om_t + (c_l \xx + \uu) \cdot \na \om & = \th_x + c_{\om} \om ,  \quad  \th_t + (c_l \xx + \uu)\cdot \na \th   =  c_{\th} \th 
\eal
\eeq
where the velocity field $\uu = (u , v)^T : \R_+^2 \times [0, T) \to \R^2_+$ is determined via the Biot-Savart law
\beq\label{eq:biot}
 - \D \phi = \om , \quad  u =  - \phi_y , \quad v  = \phi_x,
\eeq
with boundary condition $\phi(x, 0) \equiv 0$. 

An exact self-similar blowup profile solves the steady state equation of \eqref{eq:bousdy1} with \eqref{eq:biot}
\beq\label{eq:SS}
 (c_l \xx + \uu) \cdot \na \om  = \th_x + c_{\om} \om ,  \quad   (c_l \xx + \uu)\cdot \na \th   =  c_{\th} \th .
\eeq

Let $(\bar \om , \bar \th, \bar c_l , \bar c_{\om})$ be the approximate steady state 
of \eqref{eq:bousdy1} constructed in \cite{ChenHou2023a,ChenHou2023b}. 
The profile $(\bar \om, \bar \th)$ is represented using explicit basis functions $f_i(x, y)$, consisting of piecewise polynomials (B-splines), semi-analytic functions, and some explicit analytic functions; $\bar c_l , \bar c_{\om}$ are scalar self-similar exponents. The profile residual error is defined as the error in solving the equation of steady state
\beq\label{eq:bous_err}
\bga
 \bar F_{\om}   = - (\bar c_l x + \bar \uu ) \cdot \na \bar \om + \bar \th_x + \bar c_{\om} \bar\om \,  , \quad  \bar F_{\th} = - (\bar c_l x + \bar \uu ) \cdot \na \bar \th + \bar c_{\th} \bar \th , 
 \\
\olin \cF_1  \teq \bar F_{\om}, \qquad \olin \cF_2 \teq \pa_x \bar F_{\th}, \qquad \olin \cF_3 \teq \pa_y \bar F_{\th} \, .
 \ega
 \eeq

In the analysis \cite{ChenHou2023a}, we analyze the system of $(\om, \th_x, \th_y)$, which involves the error terms $\bar \cF_2$ and $ \bar \cF_3$. The term $\bar F_{\th}$ without derivative is \emph{not} used.

Due to the nonlocal Biot-Savart law, we construct the numerical approximation of the exact stream function $(-\D)^{-1} \bar \om$ by solving the Poisson equation numerically, $ - \D \bar \phi^N = \bar \om$, using the finite element method and representing $\bar \phi^N$ by explicit basis functions. 
We then construct the associated approximate velocity $ \bar \uu^N = \na^{\perp} \bar \phi^N$. Using the \emph{local} explicit basis functions, we can evaluate $\bar \om, \bar \th, 
\bar \phi^N$ and their derivatives at any point and estimate their piecewise bounds rigorously. 

The approximate stream function $\bar \phi^N$ induces the nonlocal error in solving the Poisson equation
\beq\label{eq:nloc_err}
  \bar \e = \bar \om  - (- \D) \bar \phi^N .
\eeq

We decompose the \emph{exact} velocity $\uu(\bar \om) \teq \na^{\perp}(-\D)^{-1} \bar \om$ as the approximate velocity plus error  
\beq\label{eq:vel_N}
  \uu( \bar \om) = \uu( (-\D) \bar \phi^N) + \uu( \bar \e) = \na^{\perp} \bar \phi^N + \uu(\bar \e)
  = \bar \uu^N + \uu(\bar \e) ,\quad \bar \uu^N = \na^{\perp} \bar \phi^N.
\eeq
We decompose the residual error $\bar \cF_{\bullet}$ into a local part 
involving $(\bar \om, \bar \th, \bar \phi^N)$ and a nonlocal part involving $\uu(\bar \e)$:
\footnote{
Note that the local residual error $ \bar \cF_{i}^{\lloc}$ in \eqref{def:error} 
differs from $\bar \cF_{loc, i}$ in \cite[Section 5.8]{ChenHou2023a}.
}
\beq\label{def:error}
\bal
 \bar \cF_{1}^{\lloc} & \teq 
 - (\bar c_l x + \bar \uu^N ) \cdot \na \bar \om + \bar \th_x + \bar c_{\om} \bar\om ,
  & \quad \bar \cF_1^{\nloc}  & \teq -   \uu(\bar \e)  \cdot \na \bar \om , \\
  \bar \cF_{i+1}^{\lloc}  & \teq 
  \pa_i ( - (\bar c_l x + \bar \uu^N ) \cdot \na \bar \th + \bar c_{\th} \bar \th )
  ,  & \quad 
  \bar \cF_{i+1}^{\nloc} & \teq - \pa_i   \uu(\bar \e) \cdot \na \bar \th
  -   \uu(\bar \e) \cdot \na \pa_i \bar \th , \\
\eal
\eeq
where $i=1,2$ and by a slight abuse of notation, we use $\pa_{x_1} = \pa_x, \pa_{x_2} = \pa_y$.  We refer to \cite[Section 7]{ChenHou2023a} and \cite[Appendix C,D]{ChenHou2023b} for the representations of $(\bar \om, \bar \th, \bar \phi^N)$ and the rigorous piecewise bounds.

Let $(r, \b)$ be the polar coordinates in $\R^2_+$. 
For $\max(|x|, |y|) \geq L_2, L_2 \approx 10^{15}$, 
$\bar \om, \bar \th, \bar \phi^N$ are represented by semi-analytic functions 
\beq\label{eq:decay_rate}
 \bar \om(x, y) = \bar c_1 r^{ \bar \al_1} g_1(\b),  \quad  \bar \th = \bar c_2 r^{1 + 2 \bar \al_1} g_2(\b), \quad  \bar \phi^N = \bar c_1 r^{2 + \bar \al_1} f(\b),
\eeq
 where $g_i(\b)$ and $f(\b)$ are eighth-order B--splines (see  \cite[Section 7]{ChenHou2023a}), and therefore $g_i(\b) , f(\b) \in C^{6, 1}$. 
 The profile satisfies $\bar \om, \bar \th, \bar \phi^N \in C^{4,1}$ on any compact domain,
and obeys the decay estimates 
\beq\label{eq:decay}
|\na^i \bar \om(\xx)| \les \jan{\xx}^{\bar \al_1-i}, \   |\na^i \bar \th(\xx)| \les \jan{\xx}^{1 + 2 \bar \al_1 - i}, 
\  | \na^i \bar \phi^N | \les \jan{ \xx }^{2 + \bar \al_1-i}, 
\quad \bar \al_1 \in [-0.343, -0.342] ,
\eeq
 for $i\leq 3$, where $\jan{\xx} = ( |\xx|^2 + 1)^{1/2}$.

For the local part of the error  $\bar \cF_i^{\lloc}$ in \eqref{def:error}, the decay rate $\bar \alpha_1$ in \eqref{eq:decay_rate} is chosen sufficiently close to the ratio of the scaling coefficients, $\bar c_{\om}/\bar c_l$. Consequently, the leading linear scaling terms in \eqref{def:error}, such as $\bar c_l \xx \cdot \nabla \bar \om$ and $\bar c_{\om}\bar \om$, which have the slowest decay for large $|\xx|$, 
nearly balance each other. 
\footnote{
For $\max(|x|, |y|) \geq L_2$, using \eqref{eq:decay_rate}, 
we obtain $\bar c_l \xx \cdot \na \bar \om - \bar c_{\om} \bar \om
=  \bar c_l r \pa_r \bar \om - \bar c_{\om} \bar \om
= ( \bar c_l \bar \al_1 - \bar c_{\om}) \bar \om $. Thus, in this regime, 
the local error  $\bar \cF^{\lloc}_1$ in \eqref{def:error} reduces to
$ \bar \cF^{\lloc}_1 = ( \bar c_l \bar \al_1 - \bar c_{\om}) \bar \om - \bar \uu^N \cdot \na \bar \om
+ \bar \th_x$. The first term has a very small coefficient $ \bar c_l \bar \al_1 - \bar c_{\om}$,
and the second and third term have a faster decay rate $\jan{\xx}^{2 \bar \al_1}$ 
due to \eqref{eq:decay}. The same analysis applies to 
$\bar \cF_2^{\lloc}, \bar \cF_3^{\lloc}$ in  \eqref{def:error}.
}
This cancellation is important for obtaining a small local residual error $\bar \cF_i^{\lloc}$ in the far field.

\subsection{Difference between $\bar \e$ and $\bar \cF_i$}\label{sec:differ}

For a function $g$ odd in $x$, since $\na \uu(g) \teq \na \na^{\perp}(-\D)^{-1} g$ 
is a $0$-th order singular integral operator and $\uu(g)(0) = 0$, using standard estimates for singular integral operator, decay, and regularity estimates \eqref{eq:decay}, we obtain
\beq\label{eq:err_est1}
|\na^i \bar \e | \les C_1(\bar \e) \jan{\xx}^{\bar \al_1 }, \, i\leq 1, \quad |\na \uu(\bar \e) | 
\les C_2(\bar \e)  \jan{\xx}^{  \bar \al_1 }, \quad |\xx|^{-1} |\uu(\bar \e) | 
   \les C_3(\bar \e) \jan{\xx}^{  \bar \al_1 }.
\eeq

Plugging the above estimates into \eqref{def:error} and using the decay estimates \eqref{eq:decay},
we obtain
\beq\label{eq:err_est2}
\bal
 | \bar \cF_1^{\nloc} | & \leq |\xx|^{-1} |\uu(\bar \e)| \cdot ( |\xx|  \cdot |\na \bar \om| ) 
 \leq  C_4(\bar \e) \jan{\xx}^{2 \bar \al_1}, \\
 |\bar \cF_{i+1}^{\nloc} | & \leq |\pa_{x_i} \uu( \bar \e) | \cdot |\na \bar \th|
 + |\xx|^{-1} |\uu(\bar \e)| \cdot ( |\xx|  \cdot |\na \pa_i \bar \th| )
 \leq C_{4 + i}(\bar \e)  \jan{\xx}^{3 \bar \al_1} , \quad  i = 1, 2.
 \eal
\eeq

The above constants $C_i(\bar \e)$ depend on some 
weighted $L^{\infty}, C^{1/2}$ norms of the nonlocal error $\bar \e$.

By definition, the nonlocal error $\bar \e$ and the profile residual error $\bar \cF_i$ are distinct quantities.  The nonlocal error $\bar \e$ enters the residual error $\bar \cF_i$ \eqref{def:error} 
via the nonlocal velocity $\uu(\bar \e) |\xx|^{-1}, \na \uu( \bar \e)$ 
multiplied by a  factor $ |\xx| \cdot |\na \bar \om|,
|\na \bar \th|,|\xx| \cdot | \na \pa_i \bar \th|$, which has a decaying rate  
\beq\label{eq:decay_factor}
|\xx| \cdot |\na \bar \om| \les \jan{ \xx}^{\bar \al_1},
\quad |\na \bar \th| \les \jan{ \xx}^{ 2 \bar \al_1}, \quad 
|\xx| \cdot | \na \pa_i \bar \th|  \les \jan{ \xx}^{ 2 \bar \al_1} ,
\eeq
with $\bar \al_1$ in \eqref{eq:decay} about $-\f{1}{3}$. These decaying factors are precisely what reduce the contribution of the nonlocal error to $ \bar \cF_i^{\nloc}, \bar \cF_i$ in the far field.  Since $\uu(\bar \e)$ depends nonlocally on $\bar \e$, we develop sharp functional inequalities and use appropriate  weighted norms of $\bar \e$  to bound $\uu(\bar \e), \na \uu(\bar \e)$ in \cite[Section 4]{ChenHou2023b} effectively.

\subsection{The numerical comparisons in Remarks 2.5 and 2.6 of \cite{zhang2025dimension} }

In  \cite[Remark 2.5]{zhang2025dimension}, 
the author made the following claim: 

\emph{
From the data file ``solu.w\{1,1\}'' in reference [1] of [3] in the folder file
\[
    ``Steady\_state\_pertb720\_Nlewcor4.mat''
\]
and sub folder ``solu'', there are several rectangles where the above local condition
in the sectors or condition (b') seem to hold if one uses piecewise affine interpolation
between mesh points. For instance, the rectangle given by the $x-y$ $(z^{(1)}-z^{(2)})$
mesh points $(14,719)$, $(15,719)$, $(14,718)$ and $(15,718)$. The strict maximum value
of $W = 9.02 \times 10^{-17}$ occurs at the upper right corner. Another example is the
rectangle given by the mesh points $(19,710)$, $(19,712)$, $(21,710)$, $(21,712)$. Here
interior maximum value $2.3461 \times 10^{-16}$ is reached at $(20,711)$. We comment that
in the far field, at small scales, rectangles are close to the sectorial domains.}

\emph{In addition, at mesh point $(36,720)$, we see from the \emph{solu.w\{1,1\}} file that}
\begin{equation}\label{eq:259}
\tag{2.59}
    W = 8.843970087044398 \times 10^{-17}.
\end{equation}

\emph{
From the files \emph{solu.u1\{1,2\}} and \emph{solu.u2\{2,1\}}, we find
}
\[
    \partial_{z^{(2)}} V^{(1)}
    =
    -2.850898611910781 \times 10^{-19},
    \qquad
    \partial_{z^{(1)}} V^{(2)}
    =
    1.023527319550463 \times 10^{-5}
\]
\emph{
respectively. Therefore
}
\begin{equation}\label{eq:260}
\tag{2.60}
    W
    =
    \partial_{z^{(2)}} V^{(1)}
    -
    \partial_{z^{(1)}} V^{(2)}
    \approx
    -1.023 \times 10^{-5}.
\end{equation}

\emph{
From (2.59) and (2.60) we see the difference in the order $10^{-5}$ with opposite sign.
This is beyond the accuracy of $10^{-7}$ as stated on p64 of [2] for the approximate
steady state $W$.
}

\vs{0.1in}
In \cite[Remark 2.6]{zhang2025dimension},\label{para:26} the author made the following claim:

\emph{
However, in the file solu.W\{1,1\} mentioned above, one can find
several places where $W$ is negative, including mesh point $(3,709)$, where
\[
W \approx -2.9778 \times 10^{-18}.
\]
This shows that the approximate steady state there is not an approximation
of the SSS with accuracy within $10^{-18}$ at that mesh point. Note that
(2.60), which is beyond the stated margin of error, will also lead to
a contradiction.
}

\vs{0.1in}

These interpretations are misleading because they conflate several different quantities. We clarify these distinctions in the following steps.

{\bf First}, since the velocity components \texttt{solu.u1}, \texttt{solu.u2} in the data file \texttt{solu} correspond to the grid point values of $\bar \uu^N$ defined in \eqref{eq:vel_N}, the quantity $ \partial_{z^{(2)}} V^{(1)}- \partial_{z^{(1)}} V^{(2)}$ in \eqref{eq:260} corresponds to
\[
  \pa_y \bar u_1^N -  \pa_x \bar u_2^N = - \pa_{yy} \bar \phi^N - \pa_{xx} \bar \phi^N = - \D \bar \phi^N .
\]
Thus, the difference between the value of $W$ in \eqref{eq:259} for the grid point values of $\bar \om$
and in \eqref{eq:260} for $- \D \bar \phi^N$ is the nonlocal error $\bar \e$ in \eqref{eq:nloc_err} up to a much smaller round-off error. This is a different quantity from the profile residual error $\bar \cF_i$ in \eqref{def:error}.

{\bf Second}, the quoted statement on page~64 of version~1 of
\cite{ChenHou2023a} refers to the following:

\textit{
It is extremely challenging to obtain an approximate steady state with a
sufficiently small residual error, e.g. of order $10^{-7}$, since the
solution is supported on the whole half-plane $\mathbb{R}^2_+$ and has a
slowly decaying tail in the far field, e.g.
$\omega(t,x)\sim |x|^{-1/3}$ for large $x$.
}

This sentence does not state a certified residual bound of $10^{-7}$ for the approximate steady state. Rather, it is a motivational statement explaining the difficulty of constructing an approximate steady state with a sufficiently small residual error in the presence of a slowly decaying far-field tail.

The residual error used in the rigorous analysis is the profile residual
$\bar{\mathcal F}_i$ defined in \eqref{def:error}. In particular, it is different from both the error between the approximate steady
state $\bar\omega$ and a hypothetical exact steady state, and the nonlocal error 
$\bar \e$ defined in \eqref{eq:nloc_err}.  The existence of an exact steady
state near the approximate profile is neither \emph{used} 
nor \emph{established} in the analysis of
\cite{ChenHou2023a}.

{\bf Third}, 
as explained in Section~\ref{sec:differ}, $\bar\e$ enters the profile residual error only through the
nonlocal velocity errors $|\xx|^{-1}\uu(\bar\e)$ and $\na\uu(\bar\e)$,
multiplied by additional far-field decaying factors; see \eqref{eq:err_est2}, \eqref{eq:decay_factor}.
This decay is the mechanism that makes the nonlocal contribution
$\bar\cF_i^{\nloc}$ to the profile residual much smaller in the far field than
the raw nonlocal error $\bar\e$. 
This reduction is significant at the mesh point with index $(36,720)$ (quoted above \eqref{eq:259}) used in \cite{zhang2025dimension}, which has coordinates
\footnote{
The nonuniform mesh array $\{x_i\}$ for the grid points $(x_i, x_j)$ is stored in the file \texttt{Mesh.x} within \texttt{Steady\_state\_pertb720\_Nlewcor4.mat}. The same mesh is used in both the $x$ and $y$ directions.
}
\[
(x, y) = ( x_{36}, x_{720}) \approx (0.1367, \, 2.166\cdot 10^{13} ). 
\]
At this far-field point, the additional decay factors satisfy 
\[
    \jan{ ( x_{36},  x_{720} )}^{\bar\al_1}
    \approx 2.716\cdot 10^{-5},
    \qquad
    \jan{ ( x_{36},  x_{720} ) }^{ 2 \bar\al_1}
    \approx 7.375\cdot 10^{-10}.
\]
Thus, the nonlocal error of size $10^{-5}$ quoted in \eqref{eq:260} should not be identified with either the corresponding nonlocal profile residual error or the full profile residual error.

The same decay mechanism is also exploited in the weighted $L^\infty$ energy estimates.
For $|\xx|\geq 20$, the weighted residual errors enter through terms of
the form
\footnote{
For $|\xx|\leq 20$, we estimate the residual error
$\bar\cF_i$ with correction near $\xx=0$ that improves the vanishing order of $\bar \cF_i$. We refer to
\cite[Section 5.8]{ChenHou2023a} for further discussion and do not discuss it
in this short note.
}
\footnote{
In the weighted energy estimates, we do not simply take the spatial
$L^\infty$-norms of the terms in \eqref{eq:EE}. Instead, we use the localized 
stability estimate in \cite[Lemma A.2]{ChenHou2023a}, which compares the coefficients in
the linear damping terms, the weighted residual errors, and the nonlinear terms
in order to obtain sharp stability estimates.
}
\beq\label{eq:EE}
    |\bar\cF_i \vp_i(\xx)|,\quad i\leq 3,
    \qquad
    \tau_2\mu_4 |\bar\cF_1 \vp_{g1}|,
    \qquad
    \tau_2 |\bar\cF_2 \vp_{g2}|,
    \qquad
    \tau_2 |\bar\cF_3 \vp_{g3}|.
\eeq
Here $\tau_2=0.23$, $\mu_4=0.065$, and the weights
$\vp_i,\vp_{gi}$ are chosen in \cite[Appendix C]{ChenHou2023a}. 
For each $i$, the weight $\vp_i$ decays faster than $\vp_{gi}$ for large $|\xx|$. The slowest
decay among the growing weights is
\[
    \vp_{g1}: |\xx|^{1/16},
    \qquad
    \vp_{g2}, \ \vp_{g3}: |\xx|^{1/3+10^{-8}}.
\]
These powers are chosen so that the far-field decay from 
the factors \eqref{eq:decay_factor}
is preserved in weighted estimates: 
\[
\jan{\xx}^{1/16+\bar \al_1} , \quad \jan{ \xx}^{ 1/3 +10^{-8}+2\bar\al_1} ,
\quad 
    \tfrac{1}{16}+\bar\al_1<0,
    \qquad
    \tfrac{1}{3}+10^{-8}+2\bar\al_1<0,
\]
since $\bar \al_1 \in [-0.343, -0.342] $ according to \eqref{eq:decay}.
Consequently, after the weights and parameters in \eqref{eq:EE} are included, the nonlocal contribution to the
weighted residual error still decays in the far field. Combining this decay
with the smallness of $\bar\e$, the weighted nonlocal residual errors are sufficiently
small to be treated perturbatively.

\vspace{0.05in}

{\bf Fourth}, the quoted values of the approximate profile $W$ at
the mesh points $(14,719)$, $(15,719)$, $(14,718)$, $(15,718)$ and
$(19,710)$, $(19,712)$, $(21,710)$, $(21,712)$ are all smaller than
$3\times 10^{-16}$, which are at the level of round-off error. The nonlinear
stability result established in \cite{ChenHou2023a,ChenHou2023b} controls
perturbations around the approximate profile $(\bar\omega,\nabla\bar\theta)$
only in a suitable energy norm: $E<E_*$, with $E_*=5\cdot 10^{-6}$; see 
\cite[Theorem~3]{ChenHou2023a}. Even if one could prove the existence of an \emph{exact} steady state $(\om_*, \na \th_*)$ inside this energy ball around $(\bar\omega,\nabla\bar\theta)$, 
\footnote{
For the De Gregorio and Hou-Luo model \cite{chen2019finite,chen2021HL}, 
using a time-differentiation argument, Chen-Hou-Huang established the existence of an \emph{exact} profile within the energy ball around the approximate profile for nonlinear stability estimates.
}
the bound for the perturbation $(\om_* - \bar \om, \na \th_* - \na \bar \th)$ implied by the weighted estimates $E<E_*$ is much larger than the quoted grid-point values of $\bar \om$. Therefore, the maximum of the exact steady state $(\om_*, \na \th_*)$ need not occur at comparable locations, and
these grid point values do not provide compelling support for the local sector/rectangle conditions discussed
in \cite[Remark 2.5]{zhang2025dimension}.

Regarding \hyperref[para:26]{\cite[Remark 2.6]{zhang2025dimension}}, we remark that
\cite{ChenHou2023a,ChenHou2023b} make no claim that the approximate steady state is  an approximation of an exact SSS (self-similar solution) with accuracy within $10^{-18}$ at that mesh point  ($(3, 709 )$).

\vspace{0.05in}

In summary, the comparison made in \cite{zhang2025dimension} conflates the
raw nonlocal error $\bar\e$ in \eqref{eq:nloc_err} with the profile residual
error $\bar\cF_i$ in \eqref{def:error}. The observation in \cite{zhang2025dimension} concerns the nonlocal error $\bar\e$ in the far field, whereas the relevant nonlocal profile residual
$\bar\cF_i^{\nloc}$ involves the nonlocal velocity errors
$|\xx|^{-1}\uu(\bar\e)$ and $\na\uu(\bar\e)$ multiplied by additional
far-field decay factors; see \eqref{eq:decay_factor}. These factors provide
crucial smallness and substantially reduce the contribution of the nonlocal
error to the residual in the far field. Moreover, the sentence on page~64 of
\cite{ChenHou2023a} does not assert a certified residual bound for the approximate steady state. Therefore, the comparison in \cite{zhang2025dimension} does not reveal a discrepancy with any certified residual bound. The quoted conclusion ``{\it This is beyond the accuracy of $10^{-7}$ as stated on p64 of \cite{ChenHou2023a} for the approximate steady state $W$}'' is not justified by this comparison. It rests on an incorrect identification of the raw nonlocal Poisson error with the profile residual error used in the stability analysis, and therefore gives a misleading impression of the residual estimates in \cite{ChenHou2023a,ChenHou2023b}.

\begin{remark}

In an October 2023 email exchange with Chen, Zhang pointed out that the difference between the vorticity 
profile data \texttt{solu.w\{1,1\}} and the quantity $\partial_{z^{(2)}}V^{(1)}-\partial_{z^{(1)}}V^{(2)}$, where $\partial_{z^{(2)}}V^{(1)}$ and $\partial_{z^{(1)}}V^{(2)}$ are computed from the files \texttt{solu.u1\{1,2\}} and \texttt{solu.u2\{2,1\}}, respectively, is of order $O(10^{-5})$ at mesh point $(36,720)$. He also asked whether the round-off error in this computation was around $10^{-7}$. Chen briefly responded that this difference is related to the nonlocal error in solving for the velocity, and that the velocity enters the equation only through nonlinear terms, such as $\uu\cdot\na\om$, which decay faster than the linear terms. As clarified above, this difference is the Poisson error $\bar\e$ in \eqref{eq:nloc_err}, up to a much smaller round-off error; it is distinct from the profile residual error $\bar\cF_i$ in \eqref{def:error}.

\end{remark}

\begin{remark}

The statement

\vspace{0.05in}

\emph{
This is beyond the accuracy of $10^{-7}$ as stated on p64 of [2] for the approximate steady state $W$.
}

\vspace{0.02in}

\noindent appears in Remark 2.5 of version 5 of \cite{zhang2023dimensionarXiv}. The statement

\vspace{0.05in}

\emph{
Note that (2.60), which is beyond the stated margin of error, will also lead to a contradiction.
}

\vspace{0.02in}

\noindent appears in Remark 2.6 of version 5 of \cite{zhang2023dimensionarXiv}. Both statements also appear in the published version \cite{zhang2025dimension}.

Although Zhang asked in the October 2023 correspondence whether the round-off error was around $10^{-7}$, the two affirmative conclusions quoted above were not stated in that correspondence or in the earlier email exchanges with J. Chen. These statements also did not appear in the first four arXiv versions (v1--v4) of \cite{zhang2023dimensionarXiv}. In particular, the acknowledgment in \cite{zhang2023dimensionarXiv,zhang2025dimension} of discussions with J.~Chen should not be understood as indicating his agreement with the above statements made in \cite{zhang2023dimensionarXiv,zhang2025dimension} concerning the work \cite{ChenHou2023a,ChenHou2023b}.

\end{remark}

\subsection{Figure-based evidence in  \cite[Remark 2.4]{zhang2025dimension} }

In \cite[Remark 2.4]{zhang2025dimension}, the author discusses the
applicability of \cite[Proposition 2.3(i)]{zhang2025dimension} to the
approximate profile in the Hou--Luo scenario
\cite{luo2013potentially-2,ChenHou2023a,ChenHou2023b}. Under the assumptions
of that result, asymptotically self-similar and exact self-similar profiles
are ruled out. One of the assumptions is the following:
\beq\label{eq:ass}
\text{There exist } \beta_1>0 \text{ and }
\beta_2\in(\beta_1,\pi/2) \text{ such that} 
\sup_{r>0} W(r,\beta_2)
>
\sup_{r>0} W(r,\beta_1).
\eeq
Here $(r,\beta)$ denotes the polar coordinates in $\R^2_+$.

\begin{figure}[t]
    \centering
 \begin{subfigure}{0.3\textwidth}
        \centering
        \includegraphics[width=\textwidth]{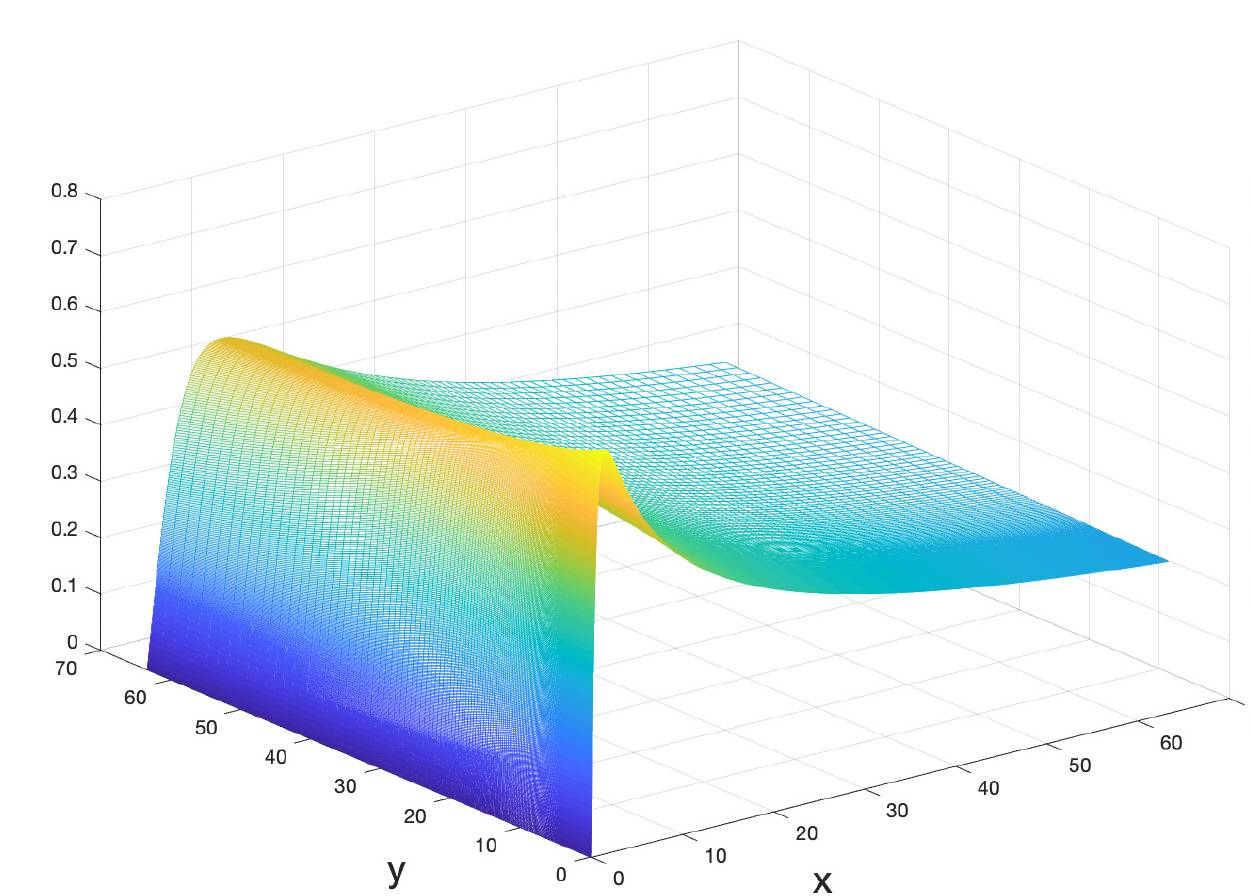}
    \end{subfigure}
    \hfill
    \begin{subfigure}{0.34\textwidth}
        \centering
        \includegraphics[width=\textwidth]{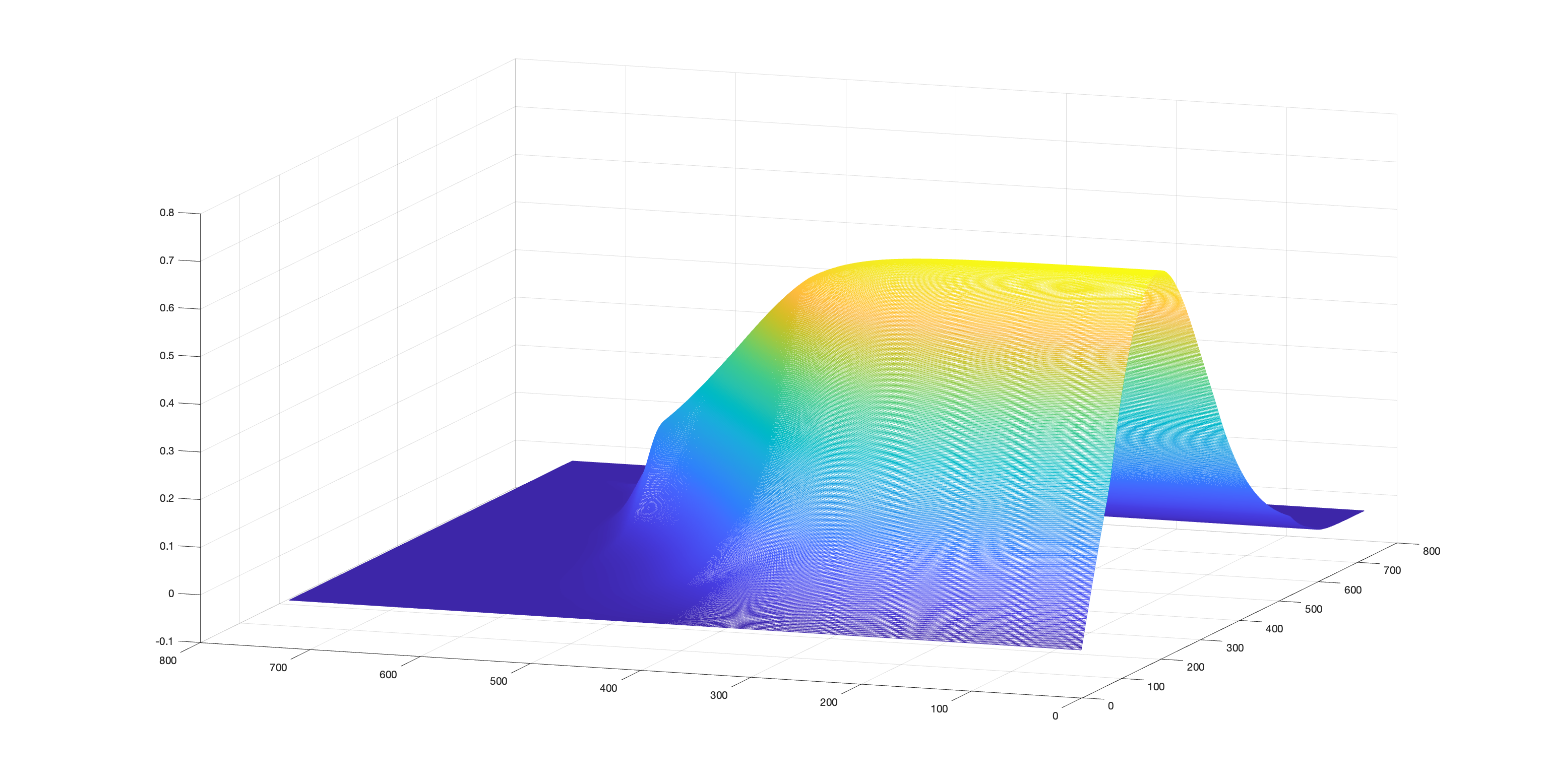}
    \end{subfigure}
    \hfill
    \begin{subfigure}{0.34\textwidth}
        \centering
        \includegraphics[width=\textwidth]{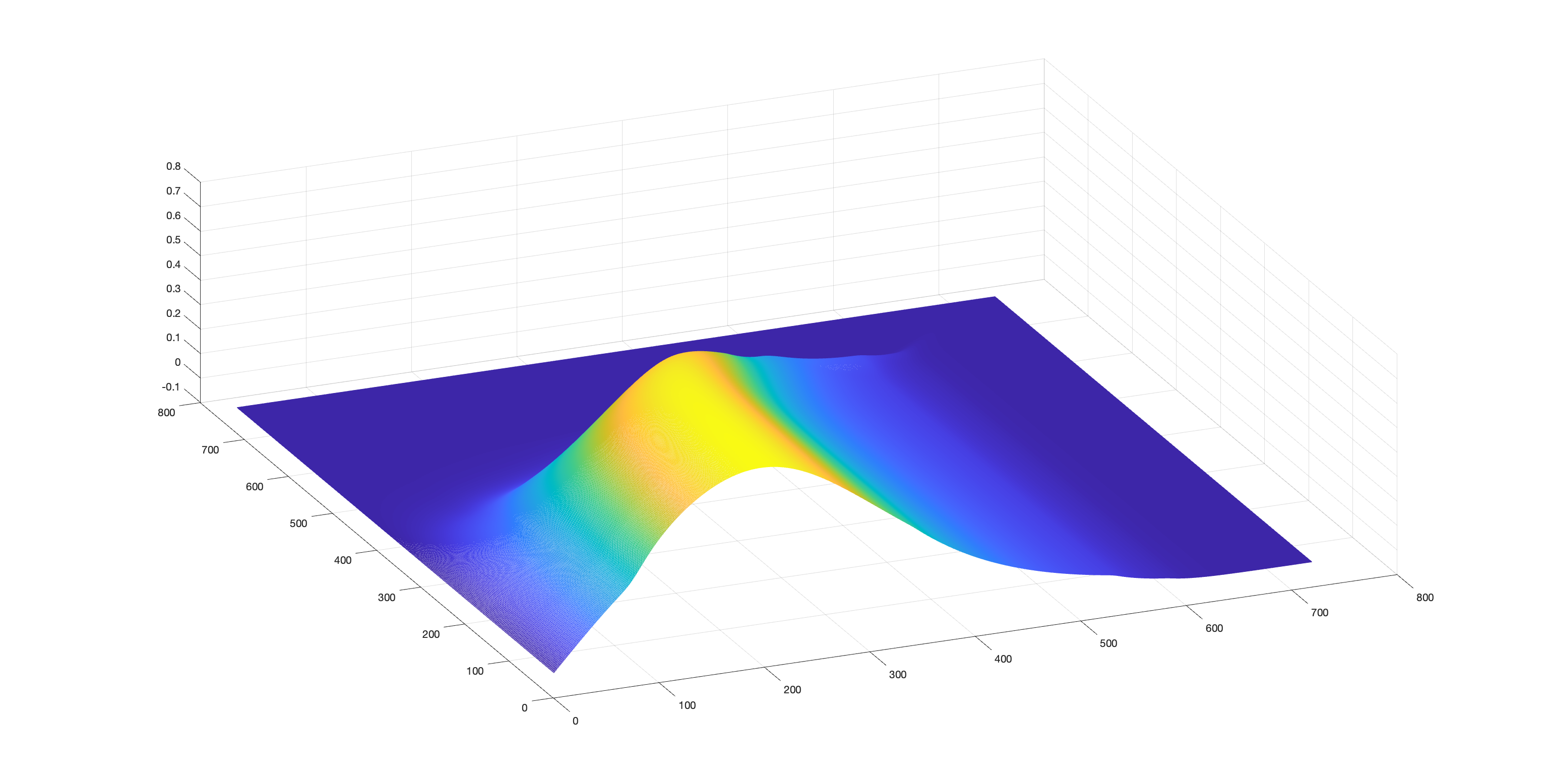}
    \end{subfigure}
    \captionsetup{width=0.94\textwidth}
    \caption{Left: reproduction of Figure~1 of version~1 of \cite{ChenHou2023a}, axes represent nonuniform mesh points $x_i$. Middle and Right: grid-point values of $\bar\omega$ plotted against mesh indices $i$ from two different perspectives. 
    Right axis corresponds to $(x,0)$ ($\beta=0$), left axis to $(0,y)$ ($\beta=\pi/2$).
}
    \label{fig:w_plots}
\end{figure}

The variable $W$ corresponds to the vorticity profile in the notation of \cite{zhang2025dimension}.   The applicability of assumption \eqref{eq:ass} to the Hou--Luo
scenario is discussed in \cite[Remark 2.4]{zhang2025dimension} based on
Figure~1 (page~9, v1) of \cite{ChenHou2023a},  reproduced in
the left panel of Figure~\ref{fig:w_plots}.
In particular, the author wrote:

\vs{0.05in}
\textit{``There is a rise of $W$ along the yellow ridge in Figure 1 after the lowest
point of the saddle, in the direction of positive $y$ axis ($z^{(2)}$ axis
here). That is why $\sup_r W(r,\theta_2)$ is larger than
$\sup_r W(r,\theta_1)$.''
}
\vs{0.01in}

Thus the numerical approximate profile $\bar\om$ is explicitly invoked as
evidence for the angular increase in \eqref{eq:ass}.\footnote{
We note that the visual interpretation of 2D plot is potentially delicate:
\cite[Remark 2.5]{zhang2025dimension} mentions the suggestion that the saddle
shape in Figure~1 of \cite{ChenHou2023a} may be caused by perspective in the
plot rather than by the underlying data.}
We therefore check whether the same grid-point data actually exhibit such an
increase for  $\bar\om$.  In the middle and right panels of Figure~\ref{fig:w_plots}, we plot the grid-point values of $\bar\om$ using the data file
\texttt{solu.w\{1,1\}} (in \texttt{Steady\_state\_pertb720\_Nlewcor4.mat}) 
from different perspectives.

Since \eqref{eq:ass} concerns the radial supremum as a function of the angle,
we test whether the data exhibit an angular increase for some pair
$\beta_1<\beta_2$. A direct verification on the numerical data would require
interpolation along many rays. As a simple sector-wise diagnostic, we instead
partition $[0,\pi/2]$ into $n=200$ equal angular sectors
$I_k \teq [\frac{\pi}{2}\frac{k-1}{n},\frac{\pi}{2}\frac{k}{n}]$,
$1\leq k\leq n$, and plot the maximal grid-point value of $\bar\om$ in each
sector:
\[
    A_k=\max\nolimits_{\beta_{ij}\in I_k}\bar\om(x_i,x_j),
    \qquad
    \beta_{ij}=\arctan(x_j/x_i),
    \qquad 1\leq k\leq n,
\]
where we take  $\beta_{ij}=\pi/2$ when $x_i=0$. \footnote{
Since $\bar \om(0, x_j)= 0$, the positive maximum is not achieved along $x=0$. 
}
Thus $A_k$ is a natural
sector-wise grid analogue of the radial supremum for 
$\bar\om$, computed from the same data \texttt{solu.w\{1,1\}} discussed above.

\begin{figure}[htbp]
    \centering
    \includegraphics[width=0.45\textwidth]{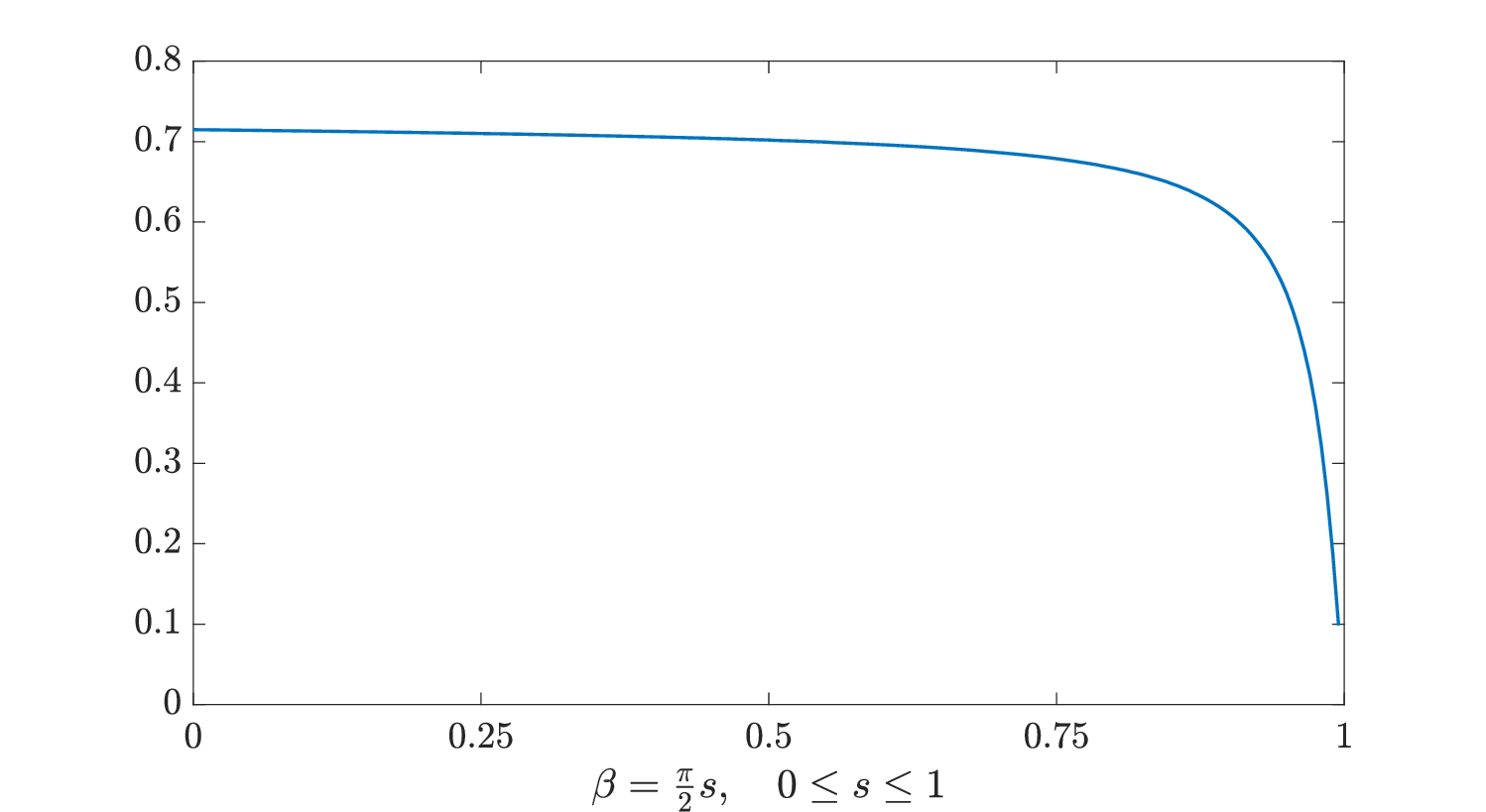}
        \captionsetup{width=0.95\textwidth}
    \caption{Sector-wise approximation $A_k$ 
    in the $k$-th angular sector $[\frac{\pi}{2}\frac{k-1}{200},\frac{\pi}{2}\frac{k}{200}]$, $1\leq k\leq 200$.  
    Larger $k$ corresponds to a larger angle.
    }
    \label{fig:plot_w4}
\end{figure}

Since the approximate profile $\bar\om$ has $C^{4,1}$ regularity and decays in
the far field, see \eqref{eq:decay} and the discussion therein, these sector-wise grid maxima provide a
meaningful numerical diagnostic for the angular profile of
$\sup_{r>0}\bar\om(r, \beta)$. With this reasonably sampled sector partition,
Figure~\ref{fig:plot_w4} shows a decreasing trend of this grid analogue as the
angle $\beta$ increases. In particular, for $\bar\om$, this diagnostic does
not reveal any angular increase of the radial maximum.

This computation tests the figure-based evidence used in
\cite[Remark 2.4]{zhang2025dimension} on the same grid-point values of
$\bar\om$. The result is that these grid-point data do not support the angular increase that \cite[Remark 2.4]{zhang2025dimension} inferred from Figure~1 of version 1 of \cite{ChenHou2023a}. Therefore, this remark does not provide a valid justification for applying the angular-increase assumption \eqref{eq:ass} to the Hou--Luo scenario, and the figure-based argument should not be used as evidence that the approximate Hou--Luo profile constructed in \cite{ChenHou2023a}  satisfies this assumption. Since a hypothetical exact  profile in the Hou--Luo scenario is not currently available, we do not attempt to verify or disprove \eqref{eq:ass} for such a profile.

\subsection{On the numerical interpretation of the one-point vanishing result in \cite[Remark 2.6]{zhang2025dimension}}

In \cite[Remark 2.6]{zhang2025dimension}, the author 
made the following claim:

\emph{
Using (2.50) and (2.61), we see that the following equations hold at the base
$z^{(2)}=0$:
}
\begin{equation}\label{eq:264}
\begin{aligned}
\left(V^{(1)}+(1-\alpha)z^{(1)}\right)\partial_{z^{(1)}}W
+W-\partial_{z^{(1)}}H^2 &=0, \\
\left(V^{(1)}+(1-\alpha)z^{(1)}\right)\partial_{z^{(1)}}H^2
+(1+\alpha)H^2 &=0.
\end{aligned}
\tag{2.64}
\end{equation}

\emph{
Under our assumption, this system has explicit solutions at the base $z^{(2)}=0$:
}
\[
\bal
H^2(z^{(1)},0)
 & =
C\exp\left(
-(1+\alpha)
\int
\left(V^{(1)}+(1-\alpha)z^{(1)}\right)^{-1}
\,dz^{(1)}
\right), \\
W(z^{(1)},0)
& =
-(1+\alpha)
\exp\left(
-\int
\left(V^{(1)}+(1-\alpha)z^{(1)}\right)^{-1}
\,dz^{(1)}
\right)
\\
\times C
\int
\exp & \left(
-\alpha
\int
\left(V^{(1)}+(1-\alpha)z^{(1)}\right)^{-1}
\,dz^{(1)}
\right) \times
\left(V^{(1)}+(1-\alpha)z^{(1)}\right)^{-2}
\,dz^{(1)}.
\eal 
\]

\emph{
From these, we see that $H^2$ can only have one zero at the origin unless it is
identically $0$, and the same holds for $W$. Note that $\alpha\approx-2$.
}

\emph{
This observation seems to indicate a potential instability in computer
calculations at the far field for the SSS. Since $W(z)$, the profile of
vorticity for the expected SSS, decays at order around $|z|^{-1/3}$ when
$|z|$ is large, a rounding error in computation may cause it to be regarded
as $0$ by the computer. According to the observation, this will force
$\omega$ and $h$ to be identically $0$ at the base.
}

\vs{0.05in}

The calculation in \cite[Remark 2.6]{zhang2025dimension} gives a rigidity property for an exact self-similar profile satisfying the boundary equations \eqref{eq:264}, under the assumptions stated there. In particular, the conclusion concerns an exact zero of an exact solution of the boundary equations.

This exact rigidity result does not, however, imply instability of the numerical construction in \cite{ChenHou2023a,ChenHou2023b}. A floating-point value that is rounded or underflowed to zero at an isolated mesh point is not an exact zero of the underlying profile. Moreover, the numerically constructed approximate profile satisfies the steady equations only up to a rigorously controlled residual error, and its restriction to the boundary may be viewed schematically as a perturbed version of \eqref{eq:264}, with additional small forcing terms arising from the residual error. The rigidity property for the exact homogeneous boundary equations need not be preserved in the presence of such forcing terms. Therefore, it cannot be applied directly to a machine-rounded grid-point value of the approximate profile.

More generally, the computer-assisted proof in \cite{ChenHou2023a,ChenHou2023b} does not require the approximate profile to agree pointwise with a hypothetical exact self-similar profile to machine precision. Round-off, interpolation, truncation, and discretization errors are incorporated into the certified residual estimates and are controlled in the weighted norms used in the stability analysis. Consequently, the one-point vanishing result for exact profiles does not establish that far-field decay 
and round-off error cause the numerical approximation to collapse to the trivial solution.

\bibliographystyle{plain}
\bibliography{selfsimilar}

\end{document}

%% file: JCmacro.tex
\usepackage{amsmath}
\usepackage{amssymb}
\usepackage{amsthm}
\usepackage{amsmath}
\usepackage{setspace}
\usepackage{xcolor}
\usepackage{fancyhdr}
\usepackage{amssymb}
\usepackage{amsthm}
\usepackage{listings}
    \usepackage{cite}
    \usepackage{tabularx}
    \usepackage{graphicx}
    \usepackage{epstopdf}    
    \usepackage{epsfig}
    \usepackage{float}
    \usepackage{ listings}
    \usepackage{appendix}
 \theoremstyle{plain}

 \newtheorem{theorem}{Theorem}[section]

 \theoremstyle{definition}

 \newtheorem{remark}[theorem]{Remark}

 \let\pa=\partial
 \let\al=\alpha
 \let\b=\beta

 \let\e=\varepsilon

 \let\f=\frac
 \let \les = \lesssim
  
 \let\om=\omega
 
 \let \th = \theta

 \let \vp = \varphi

 \let\D=\Delta

 \let \olin = \overline

 \let\teq \triangleq
 
 \let\pa=\partial
 
 \let \vs = \vspace

 \def\cF{{\mathcal F}}

 \def\na{\nabla}
 \def\la{\langle}
 \def\ra{\rangle}

 \newcommand{\bseq}{\begin{subequations}}
 \newcommand{\eseq}{\end{subequations}}

 \newcommand{\beq}{\begin{equation}}
 \newcommand{\eeq}{\end{equation}}
  \newcommand{\bal}{\begin{aligned} }
  \newcommand{\eal}{\end{aligned}}
  \newcommand{\bit}{\begin{itemize} }
  \newcommand{\eit}{\end{itemize}}
    \newcommand{\bga}{ \begin{gathered} }
  \newcommand{\ega}{ \end{gathered} }
 \newcommand{\ben}{\begin{eqnarray}}
 \newcommand{\een}{\end{eqnarray}}
 \newcommand{\beno}{\begin{eqnarray*}}
 \newcommand{\eeno}{\end{eqnarray*}}

 \newcommand{\uu}{\mathbf{u}}

 \newcommand{\xx}{\mathbf{x}}

 \newcommand{\R}{\mathbb{R}}


%% file: Euler_residual_note.bbl
\begin{thebibliography}{10}

\bibitem{chen2019finite2}
Jiajie Chen and Thomas~Y Hou.
\newblock Finite time blowup of {2D} {Boussinesq} and {3D} {Euler} equations
  with ${C}^{1,\alpha}$ velocity and boundary.
\newblock {\em Communications in Mathematical Physics}, 383(3):1559--1667,
  2021.

\bibitem{ChenHou2023a}
Jiajie Chen and Thomas~Y Hou.
\newblock Stable nearly self-similar blowup of the 2{D} {B}oussinesq and 3{D}
  {E}uler equations with smooth data {I}: {A}nalysis.
\newblock {\em arXiv preprint: arXiv:2210.07191v3 [math.AP]}, 2022.

\bibitem{ChenHou2023b}
Jiajie Chen and Thomas~Y Hou.
\newblock Stable nearly self-similar blowup of the 2{D} {B}oussinesq and 3{D}
  {E}uler equations with smooth data {II}: {R}igorous numerics.
\newblock {\em Multiscale Modeling \& Simulation}, 23(1):25--130, 2025.

\bibitem{chen2019finite}
Jiajie Chen, Thomas~Y Hou, and De~Huang.
\newblock On the finite time blowup of the {D}e {G}regorio model for the 3{D}
  {E}uler equations.
\newblock {\em Communications on Pure and Applied Mathematics},
  74(6):1282--1350, 2021.

\bibitem{chen2021HL}
Jiajie Chen, Thomas~Y Hou, and De~Huang.
\newblock Asymptotically self-similar blowup of the {H}ou--{L}uo model for the
  3{D} {E}uler equations.
\newblock {\em Annals of PDE}, 8(2):24, 2022.

\bibitem{elgindi2019finite}
Tarek~M Elgindi.
\newblock Finite-time singularity formation for ${C}^{1,\alpha}$ solutions to
  the incompressible {E}uler equations on $\mathbb{R}^3$.
\newblock {\em Annals of Mathematics}, 194(3):647--727, 2021.

\bibitem{luo2013potentially-2}
G~Luo and TY~Hou.
\newblock Toward the finite-time blowup of the 3{D} incompressible {E}uler
  equations: a numerical investigation.
\newblock {\em SIAM Multiscale Modeling and Simulation}, 12(4):1722--1776,
  2014.

\bibitem{Tristan2022}
Yongji Wang, C-Y Lai, Javier G{\'o}mez-Serrano, and Tristan Buckmaster.
\newblock Asymptotic self-similar blow-up profile for three-dimensional
  axisymmetric {E}uler equations using neural networks.
\newblock {\em Physical Review Letters}, 130(24):244002, 2023.

\bibitem{zhang2025dimension}
Qi~Zhang.
\newblock Dimension reduction of axially symmetric euler equations near maximal
  points off the axis.
\newblock {\em Transactions of the American Mathematical Society},
  378(05):3129--3156, 2025.

\bibitem{zhang2023dimensionarXiv}
Qi~S Zhang.
\newblock Dimension reduction of axially symmetric euler equations near maximal
  points off the axis.
\newblock {\em arXiv preprint arXiv:2306.09515}, 2023.

\end{thebibliography}
